\pdfoutput=1
\documentclass[a4paper,10pt,reqno]{scrartcl}
\usepackage[utf8]{inputenc}
\usepackage{latexsym,amsmath,amsfonts,amscd,amssymb,amsthm}
\usepackage{mathrsfs}
\usepackage{etoolbox}
\usepackage{pifont}
\usepackage{commath}
\usepackage{tabu}
\usepackage[numbers]{natbib}
\usepackage[resetlabels]{multibib}
\newcites{latex}{Proceedings of the conference}
\usepackage{soul}
\usepackage{eucal}
\usepackage{mathabx}
\usepackage[all]{xy}
\usepackage{graphicx,arydshln,booktabs,array}
\usepackage{tikz}
\usetikzlibrary[patterns]
\usepackage{color}
\definecolor{cobalt}{RGB}{61,89,171}
\definecolor{blue3}{RGB}{0,0,205}
\definecolor{ivory}{RGB}{238,238,224}
\usepackage[colorlinks,citecolor=blue3,linkcolor=blue3,urlcolor=blue3,pdfpagemode=UseNone]{hyperref}

\theoremstyle{plain}

\newtheorem*{theorem*}{Theorem}

\theoremstyle{definition}

\theoremstyle{remark}

\newcommand{\bC}{\mathbb{C}}

\newcommand{\bR}{\mathbb{R}}

\newcommand{\cJ}{\mathcal{J}}

\begin{document}

\title{On the history of the Hopf problem}

\author{Ilka Agricola$^*$, Giovanni Bazzoni$^\dagger$, Oliver Goertsches$^*$,\\
Panagiotis Konstantis$^*$, and S\"onke Rollenske$^*$
\\[0.5cm]
\small\texttt{\href{mailto:agricola@mathematik.uni-marburg.de}{agricola@mathematik.uni-marburg.de}},
\small\texttt{\href{mailto:gbazzoni@ucm.es}{gbazzoni@ucm.es}},
\small\texttt{\href{mailto:goertsch@mathematik.uni-marburg.de}{goertsch@mathematik.uni-marburg.de}},\\
\small\texttt{\href{mailto: pako@mathematik.uni-marburg.de}{pako@mathematik.uni-marburg.de}},
\small\texttt{\href{mailto:rollenske@mathematik.uni-marburg.de}{rollenske@mathematik.uni-marburg.de}}.
\\[0.1cm]
{\normalsize\slshape $^*$Fachbereich Mathematik und Informatik, Philipps-Universit\"at Marburg,}\\[-0.1cm]
{\normalsize\slshape Hans-Meerwein-Stra\ss e 6, 35032 Marburg}\\[-0.1cm]
{\normalsize\slshape $^\dagger$Departamento de Geometr{\'\i}a y Topolog{\'\i}a, Facultad de Ciencias Matem\'aticas,}\\[-0.1cm]
{\normalsize\slshape Universidad Complutense de Madrid, Plaza de Ciencias 3, 28040, Madrid}\\[-0.1cm]
}

\date{}
\maketitle

\begin{abstract}
This short note serves as a historical introduction to the \emph{Hopf problem}:
\begin{quote}
\hfil{}``Does there exist a complex structure on $S^6$?''\hfil{}
\end{quote}
This unsolved mathematical question was the subject of the Conference ``MAM 1 $-$ (Non-)Existence of Complex Structures on $S^6$'', which took place at Philipps-Universit\"at Marburg, Germany, between March 27\textsuperscript{th} and March 30\textsuperscript{th}, 2017.
\end{abstract}

\vskip 0.5 cm

In 1947, a few years after Hassler Whitney had proved that any differentiable manifold admits a real-analytic structure \cite{Whitney}, Heinz Hopf posed the  question whether the analogous result is true for \emph{complex structures} \cite[Page 169]{Hopf1}:
\begin{quote}
``Wir fragen nun: Kann man -- in Analogie zu diesem Sachverhalt -- jede orientierbare $M^n$ gerader Dimension (die wir als reell-analytisch annehmen d\"urfen) durch Einf\"uhrung geeigneter lokaler Koordinaten zu einer komplexen Mannigfaltigkeit machen, oder ge\-statten nur spezielle topologische Typen von Mannigfaltigkeiten die Einf\"uhrung lokaler komplex-analytischer Koordinatensysteme?''\footnote{``We ask now: can one -- in analogy to this situation -- turn any orientable $M^n$ of even dimension (which we may assume to be real-analytic), by introducing suitable local coordinates, into a complex manifold, or do only special topological types of manifolds allow the introduction of local complex-analytic coordinate systems?'' (translated by the authors)}
\end{quote}

This question is answered to the negative by Hopf himself, in the very same paper, by exhibiting infinitely many orientable even-dimensional manifolds that do not admit a complex structure, among them $S^4$ and $S^8$. In \cite[Page 170]{Hopf1} Hopf states 
\begin{quote}
``Ob die Sph\"aren $S^{2m}$ mit $m\neq 1,2,4$ komplexe Mannigfaltigkeiten sind oder nicht, habe ich nicht feststellen k\"onnen."\footnote{``I was not able to determine whether or not the spheres $S^{2m}$ with $m\neq 1,2,4$ are complex manifolds." (translated by the authors)}
\end{quote}
thus stating, implicitly, the problem for the sphere $S^6$.\\

Hopf considers the sphere bundle over a manifold $M$ whose fibre over a point $p$ consists of all directions in the tangent space at $p$, and  introduces the notion of $\cJ$-manifold: this is a manifold whose sphere bundle admits a continuous fibre-preserving self-map for which no direction is mapped to itself or its opposite. As a complex structure on a manifold $M$ induces a complex structure on each tangent space, it turns in particular $M$ into a $\cJ$-manifold. He then derives a topological obstruction to the condition of being an $\cJ$-manifold. In a footnote in \cite[Page 170]{Hopf1}, Hopf says that an alternative, but related, proof of the fact that $S^4$ is not a complex manifold was communicated to him by Charles Ehresmann; such a proof was written down a couple of years later in \cite{Ehresmann1}.\\

The notion of \emph{almost complex structure} was introduced in the same period by Ehresmann \cite[Page 3]{Ehresmann1}\footnote{These results of Ehresmann were announced in a \emph{S\'eminaire Bourbaki} in January 1947, but published only in 1949.}. The existence of an almost complex structure on a manifold $M$ immediately turns it into a $\cJ$-manifold, hence this notion is used implicitly already in Hopf's paper \cite{Hopf1}; Hopf, however, does not investigate in how far the notions of almost complex manifold (or $\cJ$-manifold) and of complex manifold are distinct. 
Apart from the fact that $S^4$ is not a complex manifold, in \cite{Ehresmann1} Ehresmann proves also a second statement which is relevant for us\footnote{Details of these proofs are contained in \cite{Ehresmann2}}, namely that a $6$-dimensional manifold with vanishing third integral homology carries an almost complex structure -- thus proving the existence of an almost complex structure on $S^6$.
We note that Ehresmann's arguments are purely topological: the existence of an almost complex structure on a manifold $M^{2n}$ amounts to a reduction of the structure group of $TM$ from $\mathrm{O}(2n)$ to $\mathrm{U}(n)$. In particular, he gives no explicit formula for the almost complex structure on $S^6$.\\

In 1947 Adrian Kirchhoff published a short paper \cite{Kirchhoff} containing two interesting results. First of all, he realized that an almost complex structure on the six-sphere can be described interpreting $S^6$ as the purely imaginary Cayley numbers of norm 1. Moreover, he proved that for $n>2$ and $n\neq 8r+6$, the sphere $S^n$ does not carry any almost complex structure. His proof is elementary and combines the following observations:  
\begin{enumerate}
\item if $S^{2n}$ is almost complex, then $S^{2n+1}$ is parallelizable;
\item as proved by Beno Eckmann \cite{Eckmann} and George W.\ Whitehead \cite{Whitehead}, if the sphere $S^n$ is parallelizable, then either $n=1$ or $n=3$, or $n$ is of the form $8r+7$.
\end{enumerate}

Of course, nowadays we know that the only parallelizable spheres are $S^1$, $S^3$ and $S^7$, but this was proved by Adams only in 1958. The last sentence of Kirchhoff's paper is
\begin{quote}
``Cependant nous ne savons pas si $S^6$ est une vari\'et\'e complexe.''\footnote{``We do not know, however, if $S^6$ is a complex manifold.'' (translated by the authors)}
\end{quote}

This is perhaps the first explicit statement of the Hopf problem. In particular, it was not clear to Kirchhoff whether the almost complex structure he had found on $S^6$ could possibly come from a honest complex structure. This question was clarified a few years later: indeed, in 1951, two papers of Eckmann and Alfred Fr\"olicher, \cite{Eckmann-Frolicher}, and Ehresmann and Paulette Libermann, \cite{Ehresmann-Libermann}, showed that the almost complex structure induced by the Cayley numbers on $S^6$ is not induced by a complex structure. Now the exceptional Lie group $\mathrm{G}_2$ is the automorphism group of the Cayley numbers, hence it is related to $S^6$ and its almost complex structure; this connection was made first in \cite{Ehresmann-Libermann}. $\textrm{G}_2$ and its relation with $S^6$ are the topic of the paper \citelatex{DF} of Cristina Draper Fontanals.\\

On the positive side, this almost complex structure on $S^6$ is part of an interesting geometric structure: the almost Hermitian structure given by this almost complex structure and the round metric is nearly K\"ahler. This type of structure was introduced by Shun-ichi Tachibana in 1959, \cite{Tachibana} and then studied systematically by Alfred Gray in 1970, \cite{Gray}; for details about the nearly K\"ahler structure on $S^6$ we refer the reader to the paper \citelatex{Agricola-Borowka-Friedrich} by Ilka Agricola, Aleksandra Bor\'owka and Thomas Friedrich.\\

In his 1949 PhD Thesis \cite{Wu1}, Wen-ts\"un Wu, using explicit relations between the Chern and Pontryagin classes which exist on a manifold endowed with an almost complex structure\footnote{Such relations, together with others involving Steenrod $p$\textsuperscript{th} powers, see \cite{Borel-Serre1,Borel-Serre2,Thom}, have been collected and developed in a paper by Friedrich Hirzebruch, \cite{Hirzebruch}.} proved that the sphere $S^{4k}$ is not an almost complex manifold, thus providing a different proof of a special case of the result of Kirchhoff. This special case of the result of Kirchhoff appears also in the influential book \emph{The topology of fibre bundles} by Norman Steenrod from 1951, \cite[Corollary 41.20]{Steenrod}. Notice that he used the terminology \emph{quasi-complex manifold} instead of almost complex manifold. In \cite[Page 209]{Steenrod}, he states:
\begin{quote}
``It seems highly unlikely that every quasi-complex manifold has a complex-analytic structure.''
\end{quote}

This question was settled 15 years later by Antonius Van de Ven: in \cite{VandeVen} he constructed the first examples of almost complex manifolds without any complex structure. His examples are 4-dimensional and it is still an open problem whether a higher-dimensional example exists. This question also appears in Shing-Tung Yau's problem list \cite{Yau}, as one of the most prominent open problems in geometry. Of course, a negative solution to the Hopf problem would yield an answer to this question.\\

Again 1951, further insights of Steenrod and John H.\ C.\ Whitehead in the problem of determining which spheres are parallelizable, together with the aforementioned ideas of Kirchhoff, show that if $n$ is not of the form $2^k-2$, then $S^n$ does not admit a complex analytic structure. Even more, it does not admit an exterior form of degree 2 which is non-singular at each point; this is precisely \cite[Corollary 1.5]{Steenrod-Whitehead}. The second statement implies that, for $n\neq 2^k-2$, $S^n$ does not admit any almost complex structure; see also \cite{Miller} for an alternative proof of this fact.\\

This is a good point to observe the following: as we recalled above, it was shown by Whitney that every differentiable manifold admits a real-analytic structure. In most of the results quoted above, the authors consider \emph{continuous} almost complex structures, hence very far from real-analytic. We should remark that the equivalence between the integrability of an almost complex structure and the vanishing of the corresponding Nijenhuis tensor, under the hypothesis that the almost complex structure is real-analytic, was known already in 1953, due to work of Eckmann and Fr\"olicher, see \cite{Eckmann2,Eckmann3,Eckmann-Frolicher}. The extension to \emph{smooth} almost complex structures, the celebrated \emph{Newlander-Nirenberg theorem}, \cite{Newlander-Nirenberg}, was proved only in 1957.\\

Announced in 1951, \cite{Borel-Serre1}, but published in 1953, \cite{Borel-Serre2}, Armand Borel and Jean-Pierre Serre proved, using Steenrod reduced $p$\textsuperscript{th} powers, that $S^2$ and $S^6$ are the only spheres admitting an almost complex structure. While the result of Hopf holds only for \emph{standard} spheres, that of Borel and Serre remains true for even-dimensional exotic spheres, i.\ e.\ manifolds of dimension $2n$ homeomorphic but not diffeomorphic to $S^{2n}$. This is an advantage of using algebraic topology. A proof of the statement that only $S^2$ and $S^6$ admit almost complex structures using K-Theory is folklore and can be found in the literature, see for instance \cite[Chapter 24]{May}. For a self-contained proof, see the paper \citelatex{Konstantis-Parton} by Panagiotis Konstantis and Maurizio Parton.\\

Thus, the question on the existence of almost complex structures on spheres was settled in the beginning of the 1950s, but the Hopf problem remained wide open. In 1955, Libermann wrote, as part of a general survey on emerging geometric structures, a first historical account on the state of the art of (almost) complex structures on spheres and the Hopf problem, see \cite{Libermann}. It is very well written and provided us with many less known bibliographical references such as \cite{Blanchard} and \cite{Wu1}. On page 204, she observes:
\begin{quote}
``Mais on ignore si la sph\`ere $S_6$ admet des structures presque complexes d\'erivant d'une structure complexe.''\footnote{``But it is not known whether the sphere $S^6$ admits almost complex structures stemming from a complex structure.'' (translated by the authors). The notation $S_n$ instead of $S^n$ for the $n$-sphere seems to have been popular, at the time, in the French literature.}
\end{quote}

After this period of intense developments, the contributions to the Hopf problem followed two main lines of research. The first concerned the non-existence of complex structures on $S^6$ satisfying some extra conditions: In 1953, Andr\'e Blanchard \cite{Blanchard} considered almost complex structures on a domain $D\subset \bR^{2n}$ which are compatible with the Euclidean metric. Using the fact that the round metric on the sphere is locally conformally flat, he proved that $S^6$ does not admit any complex structure which is compatible with the round metric. This result is usually attributed to Claude LeBrun, who rediscovered it in 1987, see \cite{LeBrun}. Together with the extension by Gil Bor and Luis Hern\'andez-Lamoneda \cite{BHL} to metrics close to the round one, this is addressed in the papers \citelatex{Ferreira,Kruglikov} by Ana Cristina Ferreira and Boris Kruglikov. A different compatibility condition, namely that with the (non-closed) fundamental form of the standard Hermitian structure on $S^6$ is the topic of the paper \citelatex{Tralle-Upmeier} by Aleksy Tralle and Markus Upmeier, which is based on an exposition of ideas of Shiing Shen Chern by Robert Bryant, see \cite{Bryant}.\\

A second line of research investigates properties of a hypothetical complex manifold $X$ with underlying differentiable manifold $S^6$, with the implicit hope to obtain either a contradiction or a hint on how to construct one. The cohomological properties of such an $X$ are the topic of the paper \citelatex{Angella} by Daniele Angella in this collection, explaining results of Gray \cite{Gray2}, Luis Ugarte \cite{Ugarte} and Andrew McHugh \cite{McHugh}. The complex analytic geometry of such an $X$, its algebraic dimension and automorphism group, was studied in \cite{CDP} by Fr\'ed\'eric Campana, Jean-Pierre Demailly and Thomas Peternell and in \cite{HKP} by Alan T.\ Huckleberry, Stefan Kebekus and Peternell; their results are presented in the paper \citelatex{LRS} by Christian Lehn, S\"onke Rollenske and Caren Schinko.\\

We end this historical introduction by mentioning two recent developments in the Hopf problem. On the one hand, there is a series of papers by Gabor Etesi, claiming the existence of a complex structure on the sphere $S^6$. Despite different proofs of this claim, based on distinct arguments, presented by Etesi, the mathematical community does not seem to have found a common opinion concerning his results - we refer to the discussions on mathoverflow \cite{mathov1,mathov2}. On the other hand, Sir Michael Atiyah posted recently a preprint on the arXiv, see \cite{Atiyah}, in which he claims to have proved the non existence of a complex structure on $S^6$. In this case as well, the community of experts does not seem to find unity, see again the mathoverflow discussions \cite{mathov3,mathov4}.


\section*{Open problems and questions}

In this section we collect questions and open problems that arose during the last talk of the conference. They reflect different ideas and hunches of the participants about the Hopf problem and possible future lines of research.

\begin{enumerate}
\item Assuming that an integrable almost complex structure $J$ exists, there exist compatible Hermitian metrics. What would ``good'' metrics be? They should belong to Grey-Hervella classe $\mathcal{W}_3\oplus\mathcal{W}_4$. Are there any of pure type?
\item What can one say about the complex geometry of a hypothetical integrable almost complex structure on $S^6$, or more generally of a compact complex 3-fold with $b_2=0$? This would entail for example existence of complex submanifolds, holomorphic fibrations, holomorphic vector bundles\ldots
\item Are there complex structures on rational homology spheres?
\item There exist infinitely many homotopy classes of almost complex structures on $\bC P^3$, indexed by the integers; the one coming from $S^6$ has ``type'' $-1$. Are any of the others integrable?
\item Does $S^6$ carry a generalized complex structure?
\item What is the K\"ahler rank\footnote{The \emph{K\"ahler rank} of a complex manifold $X$ is the maximal rank of a positive $(1,1)$ current on $X$ with analytic singularities.} of a hypothetical complex $S^6$?
\item Deformations of a hypothetical complex structure on $S^6$, moduli space of complex structures.
\item The result of Hern\'andez-Lamoneda extends Blanchard/LeBrun's result to a neighbourhood of the round metric on $S^6$. Is a similar extension of Chern's argument to a neighbourhood of the standard K\"ahler form possible?
\item Given Etesi's approach: does the existence of an integrable almost complex structure on $S^6$ have interesting physical consequences?

\end{enumerate}

\bibliographystylelatex{alpha}
\bibliographylatex{bibliography_Hopf}

\bibliographystyle{plain}
\bibliography{bibliography_Hopf}

\end{document}